\documentclass[twoside,12pt,leqno]{article}
\usepackage{amsmath, amsthm, amsfonts, amssymb, color}
\usepackage{mathrsfs}
\usepackage{fancyhdr}

\newcommand{\be}{\begin{eqnarray}}
\newcommand{\ee}{\end{eqnarray}}
\newcommand{\ce}{\begin{eqnarray*}}
\newcommand{\de}{\end{eqnarray*}}

\newtheorem{thm}{Theorem}[section]

\newtheorem{prp}[thm]{Proposition}

\theoremstyle{definition}

\definecolor{wco}{rgb}{0.5,0.2,0.3}

\numberwithin{equation}{section}
\theoremstyle{remark}
\newtheorem{rem}{Remark}[section]

\newcommand{\ua}{\uparrow}

\title{{\bf  Existence and Uniqueness of  Invariant Measures   for Stochastic
 Evolution Equations with Weakly Dissipative Drifts}
 }
\author{{\bf  Wei Liu $^a$\footnote{Corresponding author: wei.liu@uni-bielefeld.de} ,
 Jonas M. T\"{o}lle $^{b}$
}\\
{\footnotesize  $a.$  Fakult\"at f\"ur Mathematik, Universit\"at Bielefeld,
D-33501 Bielefeld, Germany}\\
  \footnotesize{$b.$ Institut f\"ur Mathematik,  Technische Universit\"{a}t Berlin, D-10623 Berlin, Germany}\\
}

\date{}
\begin{document}
\maketitle

\def\R{\mathbf R}  \def\ff{\frac} \def\ss{\sqrt} \def\BB{\mathbf
B}
\def\N{\mathbf N} \def\kk{\kappa} \def\m{{\mu}}
\def\dd{\delta} \def\DD{\Delta} \def\vv{\varepsilon} \def\rr{\rho}
\def\<{\langle} \def\>{\rangle} \def\GG{\Gamma} \def\gg{\gamma}
  \def\nn{\nabla} \def\pp{\partial} \def\tt{\tilde}
\def\d{\text{\rm{d}}} \def\bb{\beta} \def\aa{\alpha} \def\D{\mathcal D}
\def\E{\mathbb E} \def\si{\sigma} \def\ess{\text{\rm{ess}}}
\def\beg{\begin} \def\beq{\begin{equation}}  \def\F{\mathcal F}
\def\Ric{\text{\rm{Ric}}} \def\Hess{\text{\rm{Hess}}}\def\B{\mathbf B}
\def\e{\text{\rm{e}}} \def\ua{\underline a} \def\OO{\Omega}
 \def\b{\mathbf b}
\def\oo{\omega}     \def\tt{\tilde} \def\Ric{\text{\rm{Ric}}}
\def\cut{\text{\rm{cut}}}
\def\P{\mathbb P} \def\ifn{I_n(f^{\bigotimes
 n})}
\def\fff{f(x_1)\dots f(x_n)} \def\ifm{I_m(g^{\bigotimes m})}
 \def\ee{\varepsilon}
\def\pm{\pi_{{\bf m}}}   \def\p{\mathbf{p}}   \def\ml{\mathbf{L}}
 \def\C{\mathcal C}      \def\aaa{\mathbf{r}}     \def\r{r}
\def\gap{\text{\rm{gap}}} \def\prr{\pi_{{\bf m},\varrho}}
  \def\r{\mathbf r}
\def\Z{\mathbf Z} \def\vrr{\varrho} \def\ll{\lambda}
\def\L{\mathcal L}

\def\[{{\Big[}}
\def\]{{\Big]}}

\def\({{\Big(}}
\def\){{\Big)}}

\def\xy{{\|X_t-Y_t\|}}
\def\x{{\|X_t\|}}
\def\y{{\|Y_t\|}}
\def\mxy{{\mathbf{m}\left(|X_t-Y_t|^2(|X_t|\vee|Y_t|)^{r-1}\right)}}
\def\mx{{\left(\mathbf{m}\left[\left(|X_t|\vee|Y_t|\right)^{r+1}\right]\right)^{\frac{1-r}{1+r}}}}

\def\xyt{{\|X_t(x)-X_t(y)\|_H^2}}
\def\xys{{\|X_s(x)-X_s(y)\|_{H}^2}}

\def\xs{{\|X_s(x)\|_{r+1}^{1-r}}}
\def\ys{{\|X_s(y)\|_{r+1}^{1-r}}}

\def\mxys{{\mathbf{m}\left(|X_s(x)-X_s(y)|^2(|X_s(x)|\vee|X_s(y)|)^{r-1}\right)}}
\def\mxs{{\left(\mathbf{m}\left[\left(|X_s(x)|\vee|X_s(y)|\right)^{r+1}\right]\right)^{\frac{1-r}{1+r}}}}

\begin{abstract}

In this paper, a new decay estimate for a class of stochastic evolution equations with
weakly dissipative drifts is established, which directly  implies the uniqueness of invariant measures for the
corresponding transition semigroups. Moreover, the existence of invariant measures and the convergence rate
of corresponding transition semigroup to the invariant measure are  also investigated.
   As  applications, the main results  are  applied to
singular stochastic $p$-Laplace equations and
stochastic fast diffusion
equations, which solves  an open problem raised by Barbu and Da Prato in [Stoc. Proc. Appl. 120(2010), 1247-1266].
\end{abstract}



\bigbreak

\section{Introduction}

In recent years, the variational approach has been used intensively  by many authors  to analyze  
semilinear and quasilinear stochastic
partial differential equations. This approach was first
 investigated by Pardoux \cite{Pa} to study  SPDE, carried on by Krylov and
 Rozovoskii  \cite{KR} who further developed it and applied it to nonlinear filtering problems. 
 We refer to \cite{LR10,PR,Zh09} for a more detailed exposition and references.
Within this framework,  various types of both analytic and probabilistic properties  such as the large deviation principle, discretized approximation of solutions,  ergodic properties and
 existence of random attractors
have  already been established for
different  types of nonlinear SPDE (cf. \cite{GLR,L10,LR10} and references therein).  In particular, the existence and uniqueness of invariant measures and
the asymptotic behavior of the corresponding transition semigroups have been studied  for stochastic porous medium equations and stochastic
$p$-Laplace equations, see e.g. \cite{BD06,BDR,DRRW,L08a, PR,W07}.

The principal aim of this work is to show the  uniqueness of invariant measures for a class of
stochastic evolution equations with weakly dissipative drifts such as stochastic fast diffusion equations and singular
stochastic $p$-Laplace equations (where \emph{singular} means $1<p<2$ here). The existence of an invariant measure has been established by Wang
and the first named author in \cite{LW,L09a}, and recently by Barbu and Da Prato in \cite{BD10}
 for stochastic fast diffusion equations under some weaker assumptions. It is more difficult, however, to derive the uniqueness of invariant measures
for this type of stochastic equations due to the lack of strong dissipativity for the drift.
 Under some non-degeneracy assumption on the noise,  Wang and the first named author have established
the Harnack inequality for the associated transition semigroup in \cite{LW,L09a},  which implies the uniqueness of invariant measures and
some heat kernel estimates.  In this work we  prove the uniqueness of invariant measures in a more general setting and
do not assume any non-degeneracy  of the noise.
Inspired by the recent work of Es-Sarhir, von Renesse and Stannat \cite{ERS}, we establish a
decay estimate for stochastic evolution equations with weakly dissipative drifts (see $(A2)$ below), which directly implies the uniqueness of  invariant measures.
The result is applied to stochastic $p$-Laplace equations and stochastic fast diffusion equations, which also solves an open problem raised
by Barbu and Da Prato (see \cite[Remark 3.3]{BD10}).
Further applications of this new decay estimate to asymptotic behavior of corresponding transition semigroups and the construction of the corresponding
Kolmogorov operators
will be investigated in a separate paper.

Let us describe our framework in detail.
Let $H$ be a separable Hilbert space with inner product
 $\<\cdot,\cdot\>_H$ and dual
$H^*$. Let  $V$ be a reflexive  Banach space
such   that the embedding
 $V\subseteq H$ is  continuous and
dense.  Then for its dual space $V^*$ it follows that the embedding $H^*\subseteq
V^*$ is also continuous and dense.  Identifying $H$ and $H^*$ via the Riesz isomorphism we
 know that
$$V\subseteq H\equiv H^*\subseteq V^*$$
forms a so-called Gelfand triple.  If the dualization between $V^*$ and $V$ is
denoted by $_{V^*}\<\cdot,\cdot\>_V$ we have
$$ _{V^*}\<u,v\>_V= \<u,v\>_H \ \  \text{for all} \
u\in H, v\in V. $$

Suppose $\{W_t\}$ is a cylindrical Wiener process on a separable Hilbert
space $U$ w.r.t a complete filtered probability space
$(\Omega,\mathcal{F},\{\mathcal{F}_t\},\mathbb{P})$.
We consider the
following stochastic evolution  equation
 \beq
\label{1.1} d X_t=A(X_t)\,d t+B\,d W_t, \ X_0=x\in H,
 \end{equation}
 where $B$ is a Hilbert-Schmidt operator from $U$ to $H$ and
$A:  V \rightarrow V^*$
is  measurable.

 Suppose that for a fixed $\alpha>1$  there exist constants
$\delta>0$, $\beta\in (0,\alpha], \gamma\ge0$ and  $K\in\mathbb{R}$  such that the
 following
 conditions hold for all $v,v_1,v_2\in V$.
\begin{enumerate}
    \item [$(A1)$] \emph{Hemicontinuity} of $A$: The map
     $ \lambda\mapsto { }_{V^*}\<A( v_1+\lambda
 v_2),v\>_V$
     is
    continuous on $\mathbb{R}$.

    \item [$(A2)$] \emph{(Weak) dissipativity} of $A$:
$$ 2{  }_{V^*}\<A( v_1)-A( v_2), v_1-v_2\>_V
     \le -\delta \frac{\|v_1-v_2\|_H^2}{\|v_1\|_{V}^\beta+\|v_2\|_V^\beta} .$$

\item [$(A3)$] \emph{Coercivity} of $A$:
    $$ 2{ }_{V^*}\<A( v), v\>_V \le - \delta
    \|v\|_V^{\alpha} + K .$$

\item[$(A4)$] \emph{Boundedness} of $A$:
$$ \|A( v)\|_{V^*} \le  K \left( 1+
 \|v\|_V^{\alpha-1} \right) .$$
    \end{enumerate}

\begin{rem}
\begin{enumerate}
\item[(1)]
It is easy to check that  $(A1)$--$(A4)$ hold  for some concrete  examples such as the stochastic fast diffusion
equation and the singular stochastic $p$-Laplace equation (i.e. $1< p\le 2$). We refer to Section 3 for more details.

\item[(2)] $(A2)$ resembles  the following stronger dissipativity condition ($\alpha>2$):
\begin{equation}\label{strong dissipativity}
  2{  }_{V^*}\<A(v_1)-A(v_2), v_1-v_2\>_V
     \le -\delta \|v_1-v_2\|_V^\alpha, \ v_1,v_2\in V.
\end{equation}
This type of dissipativity condition holds for the stochastic porous medium equation,
 the stochastic $p$-Laplace equation ($p\ge 2$) and some other  equations with similar  degenerate drifts (cf. \cite{GLR,L08a,L10}).
 The condition $(\ref{strong dissipativity})$ has been used for the investigation of  the ergodicity in \cite{L08a},  large deviation principle  in
\cite{L10}  and the existence of random attractors in \cite{GLR}
 for a large class of stochastic evolution equations.
\end{enumerate}
\end{rem}

  Note that  $(A2)$ is stronger than the classical (weak) monotonicity condition (cf.\cite{KR,PR}),
hence for any $T>0$ and  any $x\in H$,
    $(\ref{1.1})$
    has a unique solution $\{X_t(x)\}_{t\in [0,T]}$ which is an adapted
 continuous process
    on $H$ such that $\mathbb{E}\int_0^T\|X_t\|_V^{\alpha}\,d t<\infty$
 and
$$\<X_t, v\>_H= \<x,v\>_H +\int_0^t{ }_{V^*}\<A(X_s),
    v\>_V\,d s +  \<B W_t,  v\>_H$$
holds for all $v\in V$ and
$(t, \omega)\in [0,T]\times\Omega$.

Let us define
the corresponding  transition semigroup
$$P_tF(x):= \E
F(X_t(x)),\ t\ge0, \ x\in H,$$
where $F$ is a bounded measurable function on $H$.

\begin{thm}\label{T1.0}
 Suppose $(A1)$--$(A4)$ hold for \eqref{1.1}.
\begin{enumerate}
\item[(i)] There exists  a constant $C>0$ such that
\begin{equation}\label{decay}
  \mathbb{E}  \left[  \|X_t(x)-X_t(y)\|_H^{\frac{2\alpha}{\beta}} \right]
\le  C \left(\frac{\|x-y\|_H^2}{t}\right)^{ \frac{\alpha}{\beta} }  \left(1+\frac{\|x\|_H^2}{t} +\frac{\|y\|_H^2}{t}   \right),
\ x,y\in H, t>0,
\end{equation}
where $X_t(y)$ denotes the solution of  \eqref{1.1} with starting point $y\in H$.

\item[(ii)] $\{P_t\}$ is a Feller semigroup.  Moreover, there exists $C>0$ such that
 for any Lipschitz function $F: H\rightarrow \mathbb{R}$ we have
\begin{equation}\label{Lip estimate}
|P_tF(x)-P_tF(y)|
 \le   \frac{C\mathcal{L}(F) \|x-y\|_H }{\sqrt{t} } \left(1+\frac{\|x\|_H}{\sqrt{t}} +\frac{\|y\|_H}{\sqrt{t}}   \right)^{\frac{\beta}{\alpha}},
\ x,y\in H, t>0,
\end{equation}
where $\mathcal{L}(F)$ is the Lipschitz constant of $F$.
\item[(iii)]  If $\beta\in(0, \alpha)$, then $\{P_t\}$ has at most one invariant measure.
\end{enumerate}
\end{thm}

\begin{rem} \begin{enumerate}
\item[(1)]
 This type of decay estimate  $(\ref{decay})$  is new for stochastic fast diffusion equations and singular stochastic $p$-Laplace equations.
For  stochastic porous media equations,  the following type of estimate is established in  \cite[Theorem 1.3]{DRRW} ($\alpha>2$):
$$       \|X_t(x)-X_t(y)\|_H^2\le   \|x-y\|_H^2 \wedge  \left\{ C t^{-\frac{2}{\alpha-2}} \right\}, \ t>0,\ x,y\in H.       $$
Moreover, it has been proved in \cite{L08a} that  the above estimate holds for a large class of stochastic evolution equations
with strong dissipativity condition  $(\ref{strong dissipativity})$.

\item[(2)]
 For the plane stochastic curve shortening flow (cf.\cite{ER}),  a similar type of
polynomial decay estimate is established by Es-Sarhir, von Renesse and Stannat  for the ergodic measure in \cite{ERS}.
\end{enumerate}
 \end{rem}

\begin{thm}\label{T1.1}
Suppose that $(A1)$--$(A4)$ hold with $\beta\in(0, \alpha)$  and  the embedding
$V\subseteq H$ is compact.
Then the Markov semigroup $\{P_t\}$ has an unique  invariant probability measure
 $\m$, which satisfies
 $$ \int_H\|x\|_V^\alpha \,\mu(d x)<\infty.$$
Moreover,
\begin{enumerate}
 \item[(i)]   if $\alpha\ge \sqrt{2}$, there exists $C>0$ such that  the following estimate holds:
\begin{equation}\label{ergodicity estimate}
|P_tF(x)-\mu(F)|
 \le   \frac{C\mathcal{L}(F)  \left(1+ \|x\|_H \right)   }{\sqrt{t} } \left[ 1   +\left( \frac{1+\|x\|_H}{\sqrt{t}} \right)^{\frac{\beta}{\alpha}}   \right] ,
\ x\in H, t>0;
\end{equation}
\item[(ii)] if $1<\alpha\le \sqrt{2}$, then for  any $\gamma\in \left(0, \frac{\alpha^2}{\alpha+\beta} \right]$,  there exists $C>0$ such that
\begin{equation}\label{ergodicity estimate 2}
|P_tF(x)-\mu(F)|
 \le   \frac{C |F|_\gamma \left(1+\|x\|_H^\gamma \right)}{t^{\frac{\gamma}{2}}}
 \left(  1 + \frac{1+\|x\|_H^{\frac{\beta\gamma}{\alpha}}}{t^{\frac{\beta\gamma}{2\alpha}}}  \right) ,
\ x\in H, t>0,
\end{equation}
where  $F$ is  any  $\gamma$-H\"{o}lder continuous  function and
$$|F|_\gamma:=\sup_{x\not=y\in H}\frac{|F(x)-F(y)|}{\| x-y\|_H^\gamma} .$$
\end{enumerate}
\end{thm}

The rest of the paper is organized as follows: the proofs of the main theorems are given in the next section. In Section 3, we
apply the main results to some concrete examples of SPDE.

\section{Proof of the main results}
\subsection{Proof of Theorem \ref{T1.0}}

\begin{enumerate}
\item[$(i)$]
Let $X_t(x), X_t(y)$ denote the solution of  (\ref{1.1}) starting from $x, y$ respectively.
  Then by $(A2)$ and the chain rule, we have
\ce
\|X_t(x)-X_t(y)\|_H^{2}
 &\leq& \|x-y\|_H^{2}- \delta \int_0^t
 \frac{\|X_s(x)-X_s(y)\|_{H}^2}{ \|X_s(x)\|_V^\beta +\|X_s(y)\|_V^\beta } \,d s.
\de
Note that by
$$
\frac{d}{dt}\xyt=2{  }_{V^*}\<A( X_t(x))-A( X_t(y)), X_t(x)-X_t(y)\>_V
     \le 0,
$$
the map $t\mapsto\xyt$ is decreasing.
Hence we have that
\ce
\|X_t(x)-X_t(y)\|_H^{2}
 &\leq& \|x-y\|_H^{2}- \delta \xyt  \int_0^t
 \frac{ 1}{  \|X_s(x)\|_V^\beta +\|X_s(y)\|_V^\beta }  \,d s.
\de
Furthermore,
$$  \xyt \le \|x-y\|_H^2 \left( 1+   \int_0^t
 \frac{ \delta }{  \|X_s(x)\|_V^\beta +\|X_s(y)\|_V^\beta  }\,d s    \right)^{-1} .    $$
By  Jensen's inequality,
 we get that
$$   \int_0^t
 \frac{ \delta }{   \|X_s(x)\|_V^\beta +\|X_s(y)\|_V^\beta  }  \,d s   \ge  \frac{ \delta t^2}{\int_0^t
 \left(  \|X_s(x)\|_V^\beta +\|X_s(y)\|_V^\beta  \right) \,d s  }.   $$
Which leads to
\ce
\xyt &\le& \|x-y\|_H^2 \left( 1+    \frac{ \delta t^2}{\int_0^t \left(  \|X_s(x)\|_V^\beta +\|X_s(y)\|_V^\beta \right) \,d s   }      \right)^{-1}\\
&=&  \|x-y\|_H^2 \frac{ \int_0^t \left(  \|X_s(x)\|_V^\beta +\|X_s(y)\|_V^\beta  \right) \,d s    }{  \delta t^2+ \int_0^t \left(  \|X_s(x)\|_V^\beta +\|X_s(y)\|_V^\beta \right) \,d s   }\\
&\le& \frac{ \|x-y\|_H^2  }{\delta t  }   \left(  \frac{1}{t}   \int_0^t  \|X_s(x)\|_V^\beta \,d s +
   \frac{1}{t}   \int_0^t  \|X_s(y)\|_V^\beta \,d s     \right).
\de
Using Jensen's inequality again, we get that
\begin{equation}\label{decay 1}
   \|X_t(x)-X_t(y)\|_H^{\frac{2 \alpha }{\beta}} \le   \left( \frac{ 2 \|x-y\|_H^2}{ \delta t}  \right)^{\frac{\alpha}{\beta}}
\left(  \frac{1}{t} \int_0^t \|X_s(x)\|_{V}^{\alpha}\,d s +    \frac{1}{t} \int_0^t \|X_s(y)\|_{V}^{\alpha}\,d s  \right).
\end{equation}
By applying It\^{o}'s formula to $\|\cdot\|_H^2$ and using $(A3)$,  one can easily get the estimates
$$ \mathbb{E}\left(\|X_t(x)\|_H^2+ \delta \int_0^t\|X_s(x)\|_{V}^{\alpha}  \,d s  \right)  \le \|x\|_H^2 +t \left(K+\|B\|_{HS}^2 \right); $$
$$ \mathbb{E}\left(\|X_t(y)\|_H^2+ \delta \int_0^t\|X_s(y)\|_{V}^{\alpha}  \,d s  \right)  \le \|y\|_H^2 +t \left(K+\|B\|_{HS}^2 \right),  $$
where $\|\cdot\|_{HS}$ denotes the Hilbert-Schmidt norm from $U$ to $H$.

Hence there exists  a constant $C>0$ such that
$$ \mathbb{E}  \left[  \|X_t(x)-X_t(y)\|_H^{\frac{2 \alpha }{\beta}} \right]
\le  C \left(\frac{\|x-y\|_H^2}{t}\right)^{ \frac{\alpha}{\beta} }  \left(1+\frac{\|x\|_H^2}{t} +\frac{\|y\|_H^2}{t}   \right).  $$
\item[$(ii)$] It is obvious that (\ref{decay}) implies that $\{P_t\}$ is a Feller semigroup.  Moreover, for any Lipschitz function $F: H\rightarrow \mathbb{R}$ we have
\ce
|P_tF(x)-P_tF(y)|& \le& \mathcal{L}(F)\, \mathbb{E}\,\|X_t(x)-X_t(y)\|_H \\
& \le&   \frac{C\mathcal{L}(F) \|x-y\|_H }{\sqrt{t} } \left(1+\frac{\|x\|_H}{\sqrt{t}} +\frac{\|y\|_H}{\sqrt{t}}   \right)^{\frac{\beta}{\alpha}},
\de
where $\mathcal{L}(F)$ is the Lipschitz constant of $F$ and $C>0$ is a constant (independent of $x, y, t$ and $F$).

\item[$(iii)$] Let us prove that  (\ref{Lip estimate})  is sufficient for the uniqueness of invariant measures.
It is well known that one only need to show the uniqueness of ergodic invariant measures (cf. \cite{DZ}).

In fact, if there exist two ergodic invariant measures $\mu$ and $\nu$,  then for any bounded Lipschitz function
$F$ we get in the limit  $T\to\infty$,
$$ \frac{1}{T} \int_0^T P_tF(x) \,d t\rightarrow \int_H F \,d \mu \ \text{for}  \ \mu\text{-a.e.}\ x;  $$
$$ \frac{1}{T} \int_0^T P_t F(y) \,d t\rightarrow \int_H F \,d \nu \ \text{for}  \  \nu\text{-a.e.}\  y.  $$
Since $\beta<\alpha$,   by (\ref{Lip estimate})  we have that
\ce
&  & \left|  \frac{1}{T} \int_0^T P_tF(x) \,d t-  \frac{1}{T} \int_0^T P_tF(y) \,d t    \right| \\
&\le& \frac{1}{T} \int_0^T \left|P_tF(x)-P_tF(y)\right| \,d t \\
& \le&  \frac{C \mathcal{L}(F)\|x-y\|_H }{T} \int_0^T \frac{1}{\sqrt{t}}     \left(1+\frac{\|x\|_H}{\sqrt{t}} +\frac{\|y\|_H}{\sqrt{t}}   \right)^{\frac{\beta}{\alpha}} \,d t\\
&\longrightarrow& 0 \  \text{as}\ T\rightarrow\infty.
\de
Therefore,  for any bounded  Lipschitz function $F$ on $H$ we have that
$$\int_H F \,d \mu=\int_H F\,d \nu,  $$
i.e.  $\mu=\nu$.
Therefore,   $\{P_t\}$ has at most one invariant measure.
\end{enumerate}

\subsection{Proof of Theorem \ref{T1.1}}
 Note that  $\{P_t\}$ is a Markov semigroup (cf.\cite{KR,PR}) and Feller by Theorem \ref{T1.0}.
Therefore,
 the existence  of an invariant measure $\mu$ can be proved by
 the standard Krylov-Bogoliubov procedure (cf. \cite{PR,L08a}).
 Let
 $$\mu_n:=\frac{1}{n}\int_0^n\delta_0\,P_t\,dt,\  n\geq1 ,$$
 where $\delta_0$ is the Dirac measure at $0$.

Hence for the existence of an invariant measure,  one
 only needs to verify the tightness of
 $\{\mu_n :  n\geq1\}$.

By using It\^{o}'s formula and $(A3)$, we have the following estimate:
\beq\label{estimate for invariant measure}
\x_H^2
\leq \|x\|_H^2+\int_0^t (K+\|B\|_{HS}^2-\delta\|X_s\|_{V}^{\alpha})\,ds+
2\int_0^t\langle X_s, B \,d W_s\rangle_H.
 \end{equation}

Note that $ M_t:=\int_0^t\< X_s, B\,dW_s\>_H$ is a
 martingale, then  (\ref{estimate for invariant measure})
implies that
\beq\label{estimate of mu_n}
\mu_n(\|\cdot\|_V^\alpha)=\frac{1}{n}
\int_0^n\mathbb{E}\|X_t(0)\|_{V}^{\alpha}\,dt\leq\frac{(K+\|B\|_{HS}^2)}{\delta },
\  n\geq1.
\end{equation}
Note that the embedding  $V \subseteq H$
is compact, then
 for any constant $C$ the set
$$\{x\in H: \ \|x\|_{V}\le C\}$$
is relatively compact
in $H$.  Therefore, (\ref{estimate of mu_n}) implies that
$\{\mu_n\}$ is tight, hence the limit of a convergent subsequence
provides an invariant measure $\mu$ of $\{P_t\}$.

The uniqueness of $\mu$ follows from Theorem \ref{T1.0}.  And  the concentration property for $\mu$
follows from (\ref{estimate of mu_n}), since  $\mu$ is the weak limit of $\mu_n$.

$(i)$ If $\alpha\ge \sqrt{2}$,  then  it is  also easy to show  (\ref{ergodicity estimate}) by using (\ref{Lip estimate}) and $\mu(\|\cdot\|_{V}^\alpha)<\infty$.

$(ii)$  For the case $1<\alpha\le \sqrt{2}$,   one can consider a smaller class of test function in the estimate (\ref{ergodicity estimate}). More precisely,
 for any  $\gamma$-H\"{o}lder continuous  function $F: H\rightarrow \mathbb{R}$, by H\"{o}lder's inequality and   (\ref{decay}) we have
\begin{equation*}\label{Holder estimate}
\begin{split}
 |P_tF(x)-P_tF(y)| &\le |F|_\gamma\, \E\left( \|X_t(x)-X_t(y)\|_H^\gamma  \right)\\
& \le |F|_\gamma \left[ \E\left( \|X_t(x)-X_t(y)\|_H^{\frac{2\alpha}{\beta}} \right) \right]^{\frac{\beta\gamma}{2\alpha}} \\
&\le
  \frac{C |F|_\gamma \|x-y\|_H^\gamma }{t^{\frac{\gamma}{2}} } \left(1+\frac{\|x\|_H}{\sqrt{t}} +\frac{\|y\|_H}{\sqrt{t}}   \right)^{\frac{\beta\gamma}{\alpha}},
\end{split}
\end{equation*}
where
$|F|_\gamma$
is the $\gamma$-H\"{o}lder norm of $F$ and $C>0$ is a constant (independent of $x, y, t, F$).

Hence for any $0<\gamma\le \frac{\alpha^2}{\alpha+\beta}$, we have that
\begin{equation*}
|P_tF(x)-\mu(F)|
 \le   \frac{C |F|_\gamma \left(1+\|x\|_H^\gamma \right)}{t^{\frac{\gamma}{2}}}
 \left(  1 + \frac{1+\|x\|_H^{\frac{\beta\gamma}{\alpha}}}{t^{\frac{\beta\gamma}{2\alpha}}}  \right) ,
\ x\in H, t>0,
\end{equation*}
where $F$ is  any  $\gamma$-H\"{o}lder continuous  function.

\section{Applications}
In order to verify $(A2)$ for concrete examples of stochastic evolution equations,
we first recall the following  inequality in Hilbert space proved in \cite{L09a}.
\beg{lem}\label{L3.1} Let $(H, \langle\cdot,\cdot\rangle, \|\cdot\|)$ be a
Hilbert
 space,
then for any $0<r\le 1$ we have \beq \label{6.1}
\langle\|a\|^{r-1}a-\|b\|^{r-1}b, a-b\rangle \geq r\|a-b\|^{2}
\left(\|a\|\vee \|b\|\right)^{r-1},\  a,b\in
H. \end{equation}
\end{lem}

The first example is the stochastic $p$-Laplace equation,  which arises from
geometry, plasma physics  and fluid dynamics etc (cf. \cite{Di93,L09a}).
In particular, Ladyzenskaja
suggests the $p$-Laplace equation  as a model for the
motion of non-Newtonian fluids.

Let $\Lambda$ be an open bounded domain in $\mathbb{R}^d$ with a sufficiently smooth boundary.
We consider
the following Gelfand triple
$$ V:=W^{1,p}_0(\Lambda)   \subseteq H:= L^2(\Lambda)\subseteq
 ( W^{1,p}_0(\Lambda) )^*$$
and the stochastic  $p$-Laplace equation
 \beq\label{sp}
 d X_t=\left[ \mathbf{div}(|\nabla X_t|^{p-2}\nabla X_t)
 \right]d t+B\, d W_t,\  X_0=x\in L^2(\Lambda),
  \end{equation}
 where  $p\in (1\vee \frac{2d}{2+d}, 2)$,    $B$
is a  Hilbert-Schmidt operator  on $L^2(\Lambda)$
and  $\{W_t\}$ is a cylindrical Wiener process on $L^2(\Lambda)$
w.r.t.  a complete filtered probability space
$(\Omega,\mathcal{F},\{\mathcal{F}_t\},\mathbb{P})$.

\begin{prp}
The Markov semigroup $\{P_t\}$ associated with $ (\ref{sp})$  has a unique  invariant probability measure.
\end{prp}
\begin{proof}
According to Theorem \ref{T1.1}, we need to show $(A1)$-$(A4)$ hold for (\ref{sp}).  It is well known that (\ref{sp})
satisfies
$(A1),(A3)$ and $(A4)$  with $\alpha=p$ (cf. \cite{L09a,PR}).
Let us verify $(A2)$ with $\beta=2-p$.

By Lemma \ref{L3.1} and  H\"{o}lder's  inequality   we have
\begin{equation*}
 \begin{split}
  &  {  }_{V^*}\< \mathbf{div}(|\nabla v_1|^{p-2}\nabla v_1)-
 \mathbf{div}(|\nabla v_2|^{p-2}\nabla v_2),  v_1-v_2\>_V \\
  = &  - \int_\Lambda   \left(  |\nabla  v_1|^{p-2} v_1  -   |\nabla v_2|^{p-2} v_2  \right) \left(\nabla  v_1- \nabla v_2\right) d \xi \\
\le & - (p-1) \int_\Lambda    |\nabla v_1- \nabla v_2|^2 \left(|\nabla v_1|+ |\nabla v_2|\right)^{p-2}  d \xi\\
\le & - (p-1) \frac{ \|v_1-v_2\|_{V}^2 }{ \left( \int_\Lambda (|\nabla v_1|+|\nabla v_2|)^{p}  \,d \xi  \right)^{\frac{2-p}{p}}       }\\
\le & - C \frac{\|v_1-v_2\|_H^2}{\|v_1\|_{V}^{2-p}+\|v_2\|_V^{2-p}},\ v_1,v_2\in V,
 \end{split}
\end{equation*}
where $C>0$ is some  constant derived from the  Poincar\'e   inequality used in last step.

Note that the embedding $W^{1,p}_0(\Lambda) \subseteq L^2(\Lambda)$ is compact,
hence the conclusion follows from Theorem \ref{T1.1}.  Now the proof is complete.
 \end{proof}

\begin{rem}
 In \cite{CT}, Ciotir and the second named author show  the convergence of solutions, corresponding semigroups and invariant measures for stochastic $p$-Laplace
equations as $p\to p_0$, where $p_0\in [1,2]$. In their result (see \cite[Theorem 1.5]{CT})  they  assume that the transition semigroup of  $(\ref{sp})$
 has a unique invariant measure for $p_0$.   From the above result
we know that this assumption always holds  for $p_0\in(1,2]$ (since  $d=1,2$ is assumed in \cite{CT}).
However, the uniqueness of invariant measures
in the limit  case (i.e. $p_0=1$) is still open.
\end{rem}

The second example is the stochastic fast diffusion equation, which models diffusion in plasma physics, curvature flows and
self-organized criticality in sandpile models, e.g see  \cite{BDR09,BH80,R95} and the references therein.
 We consider the following stochastic fast diffusion equation
in an open bounded domain $\Lambda$ of $\mathbb{R}^d$ with sufficiently smooth boundary (cf. \cite{BD10,LW}):
\begin{equation}\label{SFDE}
 \begin{cases}
 &d X_t(\xi) = \Delta  \left( |X_t(\xi)|^{r-1} X_t(\xi)   \right)\,d
t +B\,d W_t, \  \xi\in\Lambda,\\
 &X_t(\xi)=0, \  \forall\xi\in \partial\Lambda, \\
& X_0(\xi)=x(\xi),         \ \forall\xi\in \Lambda,
\end{cases}
\end{equation}
 where $r\in(0,1)$, $B$ is a Hilbert-Schmidt operator from $L^2(\Lambda)$ to $W^{-1,2}(\Lambda)$ and  $\{W_t\}$ is
a standard cylindrical Wiener process on $L^2(\Lambda)$ w.r.t. a complete
filtered probability space $(\OO,\F, \{\F_t\}, \P)$.

According to the classical Sobolev embedding theorem, if $r>\max\{0,\frac{d-2}{d+2}\}$, then the  embedding
$L^{r+1}(\Lambda)  \subseteq W^{-1,2}(\Lambda)$ is compact. By using
 the following Gelfand triple
$$V:=L^{r+1}(\Lambda)  \subseteq H:=W^{-1,2}(\Lambda)  \subseteq \left(L^{r+1}(\Lambda)\right)^*,$$
we can rewrite  the  stochastic fast diffusion  equation (\ref{SFDE}) into the following form:
\begin{equation}\label{FDE}
 d X_t = \Delta  \left( |X_t|^{r-1} X_t  \right)\,d
t +B \,d W_t, \ X_0=x\in H.
\end{equation}

\begin{prp} If $r\in (0 \vee \frac{d-2}{d+2}, 1)$, then
the Markov semigroup $\{P_t\}$ associated with $ (\ref{FDE})$  has a unique invariant probability measure.
\end{prp}
\begin{proof}  According to Theorem \ref{T1.1},  we only need to show $(A1)$--$(A4)$ hold for (\ref{FDE}).
 It is easy to see that  (\ref{FDE})
satisfies
$(A1),(A3)$ and $(A4)$  with $\alpha=r+1$ (cf. \cite{LW,BD10,PR}).
We will show that $(A2)$ holds with $\beta=1-r$.

Combining Lemma \ref{L3.1} with H\"{o}lder's  inequality,  we get that
\begin{equation*}
 \begin{split}
  &  {  }_{V^*}\<\Delta \left( | v_1|^{r-1} v_1  \right)-  \Delta \left( | v_2|^{r-1} v_2  \right),   v_1-v_2\>_V \\
  = &  - \int_\Lambda   \left(  | v_1|^{r-1} v_1  -   | v_2|^{r-1} v_2  \right) \left(  v_1-v_2\right) d \xi \\
\le & - r \int_\Lambda    |v_1-v_2|^2 \left(|v_1|+ |v_2|\right)^{r-1}  d \xi\\
\le & -r \frac{ \|v_1-v_2\|_{L^{r+1}}^2 }{ \left( \int_\Lambda (|v_1|+|v_2|)^{r+1} \,d \xi   \right)^{\frac{1-r}{1+r}}       }\\
\le & - C \frac{\|v_1-v_2\|_H^2}{\|v_1\|_{V}^{1-r}+\|v_2\|_V^{1-r}}, \ v_1, v_2\in L^{r+1}(\Lambda),
 \end{split}
\end{equation*}
where $C>0$ is some constant derived from  the  Sobolev  inequality used in last step.

Therefore, $P_t$ has a unique invariant measure.
\end{proof}

\begin{rem}
 For more general existence results of invariant measures for stochastic fast diffusion equations
  we refer to \cite{BD10,BDR,LW}.
If the noise in  $(\ref{FDE})$  is non-degenerate, then  the uniqueness of invariant measures, some concentration property and heat kernel estimates
has been established  in \cite{LW} by using  Harnack inequality.
 In this paper, the uniqueness of invariant measures is established without any non-degeneracy assumption on the noise,
which  answers the problem raised by Barbu and Da Prato in \cite{BD10} (see Remark 3.3 therein).
\end{rem}

\section*{Acknowledgements}
The authors would like to thank Max-K.  von Renesse for the  helpful communications.
The useful comments from the referees  are  also gratefully acknowledged. 
The first named author is supported in part by the DFG through the Internationales
 Graduiertenkolleg
``Stochastics and Real World Models'', BiBoS center and the SFB 701 ``Spectral Structures and Topological Methods in Mathematics'', Bielefeld. 
 The second named author is  supported by the DFG through the Forschergruppe 718 
``Analysis and Stochastics in Complex Physical Systems'', Berlin--Leipzig.

\end{document}